\documentclass[12pt]{amsart}

\advance\hoffset by -0.5 truecm

\usepackage{amssymb}

\newtheorem{Theorem}{Theorem}[section]
\newtheorem{Lemma}[Theorem]{Lemma}

\newtheorem{Proposition}[Theorem]{Proposition}

\def\V{\mbox{Var}}

\def\Z{{\mathbb Z}}
\def\R\re
\def\V{\bf V}

\def \re{{\mathbb R}}

\def \0{\lambda_{0}}

\def\h{{\rm h}_{\rm top}(g)}

\begin{document}
\title[Minimal entropy for 3-manifolds]{The minimal entropy problem
for 3-manifolds with zero simplicial volume}

\author[J. W. Anderson]{James W. Anderson}
\address{Faculty of Mathematical Studies\\
	University of Southampton \\
	Southampton  SO17 1BJ \\
	England} 
\email{j.w.anderson@maths.soton.ac.uk}

\author[G. P. Paternain]{Gabriel P. Paternain}
\address{Centro de Matem\'atica\\
Facultad de Ciencias\\
Igu\'a 4225\\
11400 Montevideo\\
Uruguay}
\email {gabriel@cmat.edu.uy}
\curraddr{CIMAT \\
         A.P. 402, 3600 \\
         Guanajuato. Gto. \\
         M\'exico.}
\email{paternain@cimat.mx}

\subjclass{53D25, 37D40}


\maketitle

\centerline{\it Dedicated to Jacob Palis on his sixtieth birthday}

\begin{abstract} In this note, we consider the {\em minimal entropy
problem}, namely the question of whether there exists a smooth metric
of minimal entropy, for certain classes of closed 3-manifolds.
Specifically, we prove the following two results.

\medskip
\noindent
{\bf Theorem A.} {\it Let $M$ be a closed orientable irreducible
3-manifold whose fundamental group contains a ${\Z}\oplus {\Z}$
subgroup. The following are equivalent:
\begin{enumerate}
\item the simplicial volume $\|M\|$ of $M$ is zero and
the minimal entropy problem for $M$ can be solved;
\item $M$ admits a geometric structure modelled on ${\mathbb E}^3$ or
${\rm Nil}$;
\item $M$ admits a smooth metric $g$ with $\h=0$.
\end{enumerate}}

\medskip
\noindent
{\bf Theorem B.} {\it Let $M$ be a closed orientable geometrizable
3-manifold. The following are equivalent:
\begin{enumerate}
\item the simplicial volume $\|M\|$ of $M$ is zero and
the minimal entropy problem for $M$ can be solved;
\item $M$ admits a geometric structure modelled on  ${\mathbb S}^{3}$,
${\mathbb S}^{2}\times \re$, ${\mathbb E}^3$, or ${\rm Nil}$;
\item $M$ admits a smooth metric $g$ with $\h=0$.
\end{enumerate}}
\end{abstract}


\section{Introduction and statement of results}
\label{introduction}

Let $M^{n}$ be a closed orientable $n$-dimensional manifold.  For a
smooth Riemannian metric $g$ on $M$, let ${\rm Vol}(M,g)$ denote the
volume of $M$ calculated with respect to $g$.  

Let ${\rm h}_{\rm top}(g)$ be the {\it topological entropy} of the
geodesic flow of $g$, as defined in Section \ref{topological-entropy}.
Set the {\it minimal entropy} of $M$ to be
\[{\rm h}(M):=\inf\{{\rm h}_{\rm top}(g)\: |\: g \mbox{ is a smooth
metric on $M$ with } {\rm Vol}(M,g) =1\}. \] 
A smooth metric $g_{0}$ with ${\rm Vol}(M, g_0) =1$ is {\it entropy
minimizing} if 
\[{\rm h}_{\rm top}(g_{0})={\rm h}(M). \]

The {\it minimal entropy problem} for $M$ is whether or not there
exists an entropy minimizing metric on $M$.  Say that the {\it minimal
entropy problem can be solved} for $M$ if there exists an entropy
minimizing metric on $M$. Smooth manifolds are hence naturally divided
into two classes: those for which the minimal entropy problem can be
solved and those for which it cannot.

There are a number of classes of manifolds for which the minimal
entropy problem can be solved.  For instance, the minimal entropy
problem can always be solved for a closed orientable surface $M$.  For
the 2-sphere and the 2-torus, this follows from the fact that both a
metric with constant positive curvature and a flat metric have zero
topological entropy.  For surfaces of higher genus, A. Katok \cite{Ka}
proved that each metric of constant negative curvature minimizes
topological entropy, and conversely that any metric that minimizes
topological entropy has constant negative curvature.

This result of Katok has been generalized to higher dimensions by
Besson, Courtois and Gallot \cite{BCG}, as follows.  Suppose that
$M^{n}$ ($n\geq 3$) admits a locally symmetric metric $g_{0}$ of
negative curvature, normalized so that ${\rm Vol}(M,g_0) =1$. Then
$g_{0}$ is the unique entropy minimizing metric up to isometry.
Unlike the case of a surface, a locally symmetric metric of negative
curvature on a closed orientable $n$-manifold ($n\ge 3$) is unique up
to isometry, by the rigidity theorem of Mostow \cite{mostow}.

A positive solution to the minimal entropy problem appears to single
out manifolds that have either a high degree of symmetry or a low
topological complexity.  What this means in the context of 3-manifolds
will become apparent below.  A similar phenomena is observed for
closed simply connected manifolds of dimensions 4 and 5: there are
essentially only nine manifolds for which the minimal entropy problem
can be solved and they can be explicitly listed.  These nine manifolds
share the property that their loop space homology grows polynomially
for any coefficient field, see Paternain and Petean \cite{PP}.

The goal of this note is to classify those closed orientable
geometrizable 3-manifolds with zero simplicial volume for which the
minimal entropy problem can be solved.
Specifically, in Section \ref{proof-theorem-A}, we prove:

\medskip
\noindent
{\bf Theorem A.} {\it Let $M$ be a closed orientable irreducible
3-manifold whose fundamental group contains a ${\Z}\oplus {\Z}$
subgroup. The following are equivalent:
\begin{enumerate}
\item the simplicial volume $\|M\|$ of $M$ is zero and
the minimal entropy problem for $M$ can be solved;
\item $M$ admits a geometric structure modelled on ${\mathbb E}^3$ or ${\rm
Nil}$;
\item $M$ admits a smooth metric $g$ with $\h=0$.
\end{enumerate}}

\medskip
\noindent
In Section \ref{proof-theorem-B} we prove:

\medskip
\noindent
{\bf Theorem B.} {\it Let $M$ be a closed orientable geometrizable
3-manifold. The following are equivalent:
\begin{enumerate}
\item the simplicial volume $\|M\|$ of $M$ is zero and
the minimal entropy problem for $M$ can be solved;
\item $M$ admits a geometric structure modelled on  ${\mathbb S}^{3}$, ${\mathbb
S}^{2}\times \re$, ${\mathbb E}^3$, or ${\rm Nil}$;
\item $M$ admits a smooth metric $g$ with $\h=0$.
\end{enumerate}}


\medskip

Recall that the {\it simplicial volume} of a closed orientable
manifold $M$ is defined as the infimum of $\sum_{i}|r_{i}|$ where the
$r_{i}$ are the coefficients of a real cycle that represents the
fundamental class of $M$.  For 3-manifolds, the positivity of the
simplicial volume (which is a homotopy invariant) is closely related
to the existence of compact hyperbolizable submanifolds in $M$.  This
is described in more detail in Section \ref{simplicial-volume}. 

We close the introduction by describing some of the elements of the
proofs of Theorems A and B, and by describing a conjectural picture.
We will see in Section \ref{preliminaries} that a closed orientable
geometrizable 3-manifold $M$ has zero simplicial volume if and only if
$M$ has zero minimal entropy.  Therefore, the minimal entropy problem
can be solved if and only if $M$ admits a smooth metric with zero
topological entropy. This is in turn forces the fundamental group of
$M$ to have subexponential growth.  We end up showing that this can
occur only if $M$ admits one of the four geometric structures listed
in the statement of Theorem B.  On the other hand, it is a
calculation that the manifolds in the statement of Theorem B carry a
metric of zero entropy.  The proof of Theorem A follows a similar
line, and makes use of the remarkable theorem, due essentially to
Thurston, that a manifold satisfying the hypotheses of the theorem is
geometrizable. The precise definition of geometrizable manifold is given
in Subsection 2.4. Thurston's geometrization conjecture states that
every closed orientable 3-manifold is geometrizable.

From this discussion and the above mentioned result of Besson,
Courtois and Gallot it seems quite reasonable to speculate that
the following statement holds: 

\medskip
\noindent
{\it Let $M$ be a closed orientable geometrizable 3-manifold.  Then,
the minimal entropy problem for $M$ can be solved if and only if $M$
admits a geometric structure modelled on ${\mathbb S}^{3}$, ${\mathbb
S}^{2}\times \re$, ${\mathbb E}^3$, ${\rm Nil}$, or ${\mathbb
H}^{3}$}.

\medskip
\noindent
Indeed, suppose that the simplicial volume of $M$ were {\it not}
zero. This would imply that $M$ contains a disjoint collection
$H_1,\ldots, H_p$ of compact submanifolds whose interiors each admit a
complete hyperbolic structure of finite volume.  In particular, it
should be that the minimal entropy of $M$ is the maximum of the
minimal entropies of the $H_k$. It would then seem reasonable that an
entropy minimizing metric on $M$ would try to be as hyperbolic as
possible on the interiors of the $H_k$ and would try to as much one of
the other seven standard 3-dimensional geometries as possible on the
components of $M- (H_1\cup\cdots\cup H_p)$.  However, it would seem
that the minimizer would have to be singular along the $\partial H_k$,
and so there should be no metric of minimal entropy.  Unfortunately,
we do not yet know how to make this argument rigorous.  

\section{Preliminaries}
\label{preliminaries}

The purpose of this Section is to present some of the basic material
from 3-manifold theory that we will need.  We refer the interested
reader to Hempel \cite{hempel} for a more detailed introduction to
3-manifold topology.  For a more detailed description of Seifert
fibered spaces, and of the torus decomposition and the geometrization
of 3-manifolds, we refer the interested reader to the survey articles
of Scott \cite{S} and Bonahon \cite{bonahon-chapter}, and the
references contained therein.

\subsection{3-manifold basics}
\label{prelim-3mfld}

We begin with some basic definitions.  A 3-manifold is {\it closed} it
it is compact with empty boundary.

An embedded 2-sphere ${\mathbb S}^2$ in a 3-manifold $M$ is {\it
essential} if $M$ does not bound a closed 3-ball in $M$.  A 3-manifold
is {\it irreducible} if it contains no essential 2-sphere.  

A 3-manifold is {\it prime} if it cannot be decomposed as a
non-trivial connected sum.  That is, $M$ is prime if for every
decomposition $M =M_1 \# M_2$ of $M$ as a connected sum, one of $M_1$
or $M_2$ is homeomorphic to the standard 3-sphere ${\mathbb
S}^3$.  Every irreducible 3-manifold is prime, but the converse does
not hold.  However, the only closed orientable 3-manifold that is
prime but not irreducible is ${\mathbb S}^2\times {\mathbb S}^1$.

We note here that if the closed orientable 3-manifold $M$ contains a
non-separating essential 2-sphere, then $M$ can be expressed as the
connected sum $M =P \# ({\mathbb S}^2\times {\mathbb S}^1)$ for some
3-manifold $P$.  Hence, in what follows, we need only consider
separating essential 2-spheres in 3-manifolds.

There is an inverse to the operation of connected sum for 3-manifolds,
called the {\it prime decomposition}.  The following statement is
adapted from Bonahon \cite{bonahon-chapter}, and follows from work of
Kneser \cite{Kn} and Milnor \cite{Mil2}.

Let $M$ be a closed orientable 3-manifold.  Then, there exists a
compact 2-submanifold $\Sigma$ of $M$, unique up to isotopy, so that
two conditions hold.  First, each component of $\Sigma$ is an embedded
essential separating 2-sphere, and the 2-spheres in $\Sigma$ are
pairwise non-parallel, in that no two 2-spheres in $\Sigma$ bound an
embedded ${\mathbb S}^2\times [0,1]$ in $M$.  Second, if $M_0,
M_1,\ldots, M_q$ are the closures of the components of $M -\Sigma$,
then $M_0$ is homeomorphic to the 3-sphere ${\mathbb S}^3$ minus
finitely many disjoint open 3-balls; while for $k\ge 1$, each $M_k$
contains a unique component of $\Sigma$, and every separating
essential 2-sphere in $M_k$ is parallel to $\partial M_k$.  

The {\it prime decomposition} of $M$ is the collection of 3-manifolds
that results by taking the complements of the 2-submanifold $\Sigma$
in $M$ as just described, and filling in each 2-sphere boundary
component of $M_0, M_1,\ldots, M_p$ with a 3-ball.  Each of the
resulting 3-manifolds is then prime.  (Note that both ${\mathbb S}^3$
and ${\mathbb S}^2\times {\mathbb S}^1$ have trivial prime
decompositions, as they do not contain a separating essential
2-sphere.)  The prime decomposition is one of two standard
decompositions of a closed orientable 3-manifold.

In general, a closed orientable embedded surface $S$ in a 3-manifold
$M$ is {\it 2-sided} if there exists an embedding $f$ of $S\times
[-1,1]$ into $M$ so that $f(S\times \{ 0\}) =S$.  A closed orientable
embedded surface $S$ in a 3-manifold $M$ is {\it incompressible} if
the fundamental group of $S$ is infinite and if the inclusion
$S\hookrightarrow M$ induces an injection on fundamental groups.  An
incompressible surface $S$ is {\it essential} if $S$ is not homotopic
into $\partial M$. 

A compact orientable irreducible 3-manifold $M$ is {\it sufficiently
large} if it contains a 2-sided incompressible surface.  Sufficiently
large 3-manifolds are also known as {\it Haken} 3-manifolds.  

\subsection{Seifert fibered spaces}
\label{seifert-definition}

A {\it Seifert fibration} of a 3-manifold $M$ is a decomposition of
$M$ into disjoint simple closed curves, called the {\it fibers} of the
fibration, so that each fiber $c$ has a neighborhood $U$ in $M$ of
the following form:  $U$ is diffeomorphic to the quotient of ${\mathbb
S}^1\times {\mathbb B}^2$ by the free action of a finite group action
respecting the product structure, where the fibers of the fibration
correspond to the curves $\{ x\}\times {\mathbb B}^2$ for $x\in
{\mathbb S}^1$.  (In this note, we only consider Seifert fibrations of
closed 3-manifolds and of 3-manifolds without boundary that are
homeomorphic to the interior of a compact 3-manifold with 2-torus
boundary components.)

Since we are considering only orientable 3-manifolds in this note, the
quotient of ${\mathbb S}^1\times {\mathbb B}^2$ in the above
definition can be obtained from $[0,1]\times {\mathbb B}^2$ by gluing
$(0,z)$ to $(1, z^{q/p})$, where $p$ and $q$ are relatively prime
integers.  A fiber is an {\it regular fiber} if it has a neighborhood
diffeomorphic to ${\mathbb S}^1\times {\mathbb B}^2$, and is a {\it
singular fiber} otherwise.  Note that the singular fibers of a Seifert
fibration are necessarily isolated.

Let $S$ be the space of fibers of a Seifert fibration of a 3-manifold
$M$, equipped with the quotient topology coming from the projection
map $p: M\rightarrow S$.  We often refer to $S$ as the {\it base
orbifold} of the Seifert fibered space $M$.  Using the neighborhoods
of the fibers in $M$, we see that $S$ is an orientable surface with
one cone point for each singular fiber.  

Let $p_1,\ldots, p_s$ be the cone points on $S$, and let $n_j$ be the
{\it order} at the cone point $p_j$, so that a neighbhorhood of $p_j$
is diffeomorphic to the quotient of ${\mathbb B}^2/ {\mathbb
Z}_{n_j}$, where ${\mathbb Z}_{n_j}$ acts by rotation.  The {\it
orbifold Euler characteristic} $\chi(S)$ of $S$ is the quantity 
\[ \chi(S) =2 -2\: {\rm genus}(S) - \sum_{k=1}^s \left( 1
-\frac{1}{n_j} \right). \]

(This formula is also valid in the case that $M$ is a 3-manifold
without boundary that is homeomorphic to the interior of a compact
3-manifold with 2-torus boundary components.  In this case, the base
orbifold has punctures as well as cone points, and we view each
puncture as a cone point of infinite order.)

There are two cases of particular interest.  In the case that $\chi(S)
<0$, $S$ has a hyperbolic structure, so that we can express $S$ as the
quotient $S ={\mathbb H}^2/\Gamma$, where ${\mathbb H}^2$ is the
hyperbolic plane and $\Gamma$ is a discrete subgroup of ${\rm
Isom}({\mathbb H}^2)$, where the fixed points of the action of
$\Gamma$ descend to the cone points on $S$.  We refer to $\Gamma$ as
the {\it orbifold fundamental group} of $S$.  In this case, we have
that $\Gamma$ contains a free subgroup of rank 2, and in particular
$\Gamma$ contains an element of infinite order. 

In the case that $\chi(S) =0$, $S$ has a Euclidean structure, so that
we can express $S$ as the quotient $S ={\mathbb E}^2/\Gamma$, where
${\mathbb E}^2$ is the Euclidean plane and $\Gamma$ is a discrete
subgroup of ${\rm Isom}({\mathbb E}^2)$, where the fixed points of the
action of $\Gamma$ descend to the cone points on $S$.  As above, we
refer to $\Gamma$ as the {\it orbifold fundamental group} of $S$.   In
this case, we have that $\Gamma$ contains an element of infinite
order, but not a non-trivial free subgroup.



In both of these cases, the orbifold fundamental group of the base
orbifold $S$ of the Seifert fibered space $M$ is a subgroup of
$\pi_1(M)$.  In fact, there is a short exact sequence
\[ 1\rightarrow {\mathbb Z} \rightarrow \pi_1(M) \rightarrow \pi_1(S)
\rightarrow 1, \]
where $\pi_1(S)$ is the orbifold fundamental group of $S$ and where
${\mathbb Z}$ is generated by any regular fiber of the Seifert
fibration. 

The following follows immediately from this discussion.

\begin{Lemma} Let $M$ be a Seifert fibered space as above with base
orbifold $S$.  If $\chi(S) \le 0$, then $\pi_1(M)$ contains a
${\mathbb Z}\oplus {\mathbb Z}$ subgroup.
\label{seifert-fundamental-group}
\end{Lemma}

\begin{proof} The proof of Lemma \ref{seifert-fundamental-group} is
standard, but we sketch it here for the sake of completeness.  Let $p:
M\rightarrow S$ be the quotient map.  Since $\chi(S)\le 0$, there is
a closed curve $c$  curve, not necessarily simple, on $S$ that
represents an infinite order element of the orbifold fundamental group
of $S$.  Let $T =p^{-1}(c)$ in $M$ be the subset of $M$ that consists
of all the fibers in $M$ corresponding to points of $c$. Then, $T$ is
an incompressible 2-torus in $M$, though not necessarily embedded.
However, this is sufficient to guarantee that there exists a ${\mathbb
Z}\oplus {\mathbb Z}$ subgroup of $\pi_1(M)$, namely the fundamental
group of $T$.
\end{proof}

\subsection{ The torus decomposition}
\label{section-torus-decomposition}

Let $M$ be a closed orientable irreducible 3-manifold with infinite
fundamental group.  There is then a canonical decomposition of $M$
along embedded essential 2-tori, due to Jaco and Shalen
\cite{jaco-shalen} and Johannson \cite{johannson}.  (Note that the
restriction to irreducible 3-manifolds causes no loss of generality, 
as we may first apply the prime decomposition to $M$, as described in
Section \ref{prelim-3mfld}.  Also, we tend to not take the torus
decomposition of ${\mathbb S}^2\times {\mathbb S}^1$.)  The statement
given below is adapted from Theorem 3.4 of Bonahon
\cite{bonahon-chapter}. 

\begin{Theorem} \cite{bonahon-chapter} Let $M$ be a closed
orientable irreducible 3-manifold.  Then, up to isotopy, there is a
unique compact 2-submanifold $T$ of $M$ such that:
\begin{enumerate} 
\item every component of $T$ is a 2-sided essential 2-torus; 
\item every component of $M -T$ either contains no essential embedded
2-torus or Klein bottle, or else admits a Seifert fibration (or
possibly both);
\item property (2) fails when any component of $T$ is removed.
\end{enumerate}
\label{torus-decomposition}
\end{Theorem}

We refer to this 2-submanifold $T$ as the {\it torus decomposition} of
$M$.

There are several things to note.  Condition $(3)$ implies that no two
of the 2-tori in the torus decomposition are isotopic.  Moreover,
every ${\mathbb Z}\oplus {\mathbb Z}$ subgroup of $\pi_1(M)$ is
conjugate into the fundamental group of some component of $T$. 

Let $M$ be a compact orientable 3-manifold, and let $M_0, M_1,\ldots,
M_p$ be the components of its prime decomposition.  Let $T_k$ be the
torus decomposition of $M_k$.  Say that $M$ is a {\it graph manifold}
if, for each $1\le k\le p$, every component of $M_k -T_k$ admits a 
Seifert fibration .  Clearly, every Seifert fibered space is trivially
a graph manifold.  Also, every 2-torus bundle over ${\mathbb S}^1$ is
a graph manifold.

Theorem \ref{torus-decomposition} is a small part of the machinary of
the {\it characteristic submanifold} of a 3-manifold developed by Jaco
and Shalen and by Johannson.  Note that this discussion includes the
possibility that the torus decomposition $T$ is empty, even though
$\pi_1(M)$ may contain a ${\mathbb Z}\oplus {\mathbb Z}$ subgroup.

A closely related result is the following torus theorem.  For a
discussion and proof of this result, see Scott \cite{scott-review}.

\begin{Theorem}\cite{scott-review} Let $M$ be a closed orientable
irreducible 3-manifold whose fundamental group contains a ${\mathbb
Z}\oplus {\mathbb Z}$ subgroup.  Then, either $M$ contains an
incompressible embedded 2-torus or $M$ is a Seifert fibered space.
\label{link}
\end{Theorem}


\subsection{Geometric structures and geometrization}
\label{section-geometric-structures}

A {\it 3-dimensional geometry} is a pair $(X,G)$, where $X$ is a
simply connected Riemannian 3-manifold with a complete homogeneous
metric and $G$ is a maximal transitive group of orientation-preserving
isometries of $X$, with the proviso that there exists a subgroup $H$
of $G$ with compact quotient $X/H$.  Note that since $G$ is a maximal
group of isometries, it suffices to specify $X$ and set $G ={\rm
Isom}(X)$.

It is a result of Thurston that there exist exactly eight
3-dimensional geometries, namely ${\mathbb E}^3$, ${\mathbb S}^3$,
${\mathbb H}^3$, ${\mathbb S}^2\times {\mathbb R}$, ${\mathbb
H}^2\times {\mathbb R}$, $\widetilde{{\rm SL}_2}$, ${\rm Nil}$, and
${\rm Sol}$, with their respective groups of (orientation preserving)
isometries.  (A proof of this result, and a detailed description of
the eight geometries, is given in Scott \cite{S}.) 

Let $M$ be an orientable 3-manifold that is homeomorphic to the
interior of a compact 3-manifold with 2-torus boundary components.
(This includes the possibility that $M$ is closed.)  Say that $M$ {\it
admits a geometric structure modelled on $X$} if $M$ is diffeomorphic
to the quotient $X/\Gamma$, where $X$ is one of the eight
3-dimensional geometries and $\Gamma$ is a fixed point free subgroup
of ${\rm Isom}(X)$.  It is known that if a 3-manifold admits a
geometric structure, then it admits a unique geometric structure.

More generally, let $M$ be a closed orientable irreducible
3-manifold with torus decomposition $T$.  Say that $M$ is {\it
geometrizable} if each component of $M -T$ admits a geometric
structure.  (Note that we do not require that different components
of $M -T$ admit the same geometric structure.)  

Finally, say that a closed orientable 3-manifold is {\it
geometrizable} if every component of its prime decomposition is
geometrizable.  (This causes no difficulties, as ${\mathbb S}^2\times
{\mathbb S}^1$, which may arise as a component of the prime
decomposition but is not irreducible, admits a geometric structure
modelled on ${\mathbb S}^2\times {\mathbb R}$.)

Thurston's {\it geometrization conjecture} states that every closed
orientable 3-manifold is geometrizable.  For a more complete
discussion of the geometrization conjecture, see Scott \cite{S},
Bonahon \cite{bonahon-chapter}, or Thurston \cite{thurston-bulletin}. 

There are a number of manifolds for which the geometrization
conjecture is known to be true.  If $M$ is a closed orientable
irreducible sufficiently large 3-manifold, then $M$ is geometrizable;
this is Thurston's geometrization theorem; see Morgan \cite{morgan} or
Otal \cite{otal} for a discussion of this theorem.  

In particular, if $M$ has a non-empty torus decomposition, then it is
geometrizable.  In this case, each component of the complement of the
torus decomposition of $M$ either is a Seifert fibered space or admits
a hyperbolic structure, that is the geometric structure modelled on
${\mathbb H}^3$.  We encode in the following theorem the parts of this
discussion we make the most use of.

\begin{Theorem} Let $M$ be a closed orientable irreducible
sufficiently large 3-manifold.  Then, $M$ admits a torus decomposition
$T$.  Moreover, each component of $M-T$ either is a Seifert fibered
space or admits a hyperbolic structure.
\label{hypstr}
\end{Theorem}

Additionally, the geometrization of Seifert fibered spaces, and in
fact of irreducible graph manifolds, is completely understood.

\begin{Theorem} \cite[Theorem 5.3]{S} Let $M$ be a closed orientable
3-manifold.  Then,
\begin{enumerate}
\item $M$ possesses a geometric structure modelled on ${\rm Sol}$ if
and only if $M$ is finitely covered by a 2-torus bundle over ${\mathbb
S}^{1}$ with hyperbolic glueing map;
\item $M$ possesses a geometric structure modelled on one of ${\mathbb
S}^3$, ${\mathbb E}^3$, ${\mathbb S}^2\times {\re}$, ${\mathbb
H}^{2}\times\re$, $\widetilde{{\rm SL}}_{2}$, or ${\rm Nil}$ if and
only if $M$ is a Seifert fibered space.
\end{enumerate}
\label{scott1}
\end{Theorem}

We note here that the two unresolved cases of the geometrization
conjecture are that the fundamental group of $M$ is finite, in which
case $M$ should admit a geometric structure modelled on ${\mathbb
S}^3$ [the Poincar\'e conjecture and the spherical space form
problem], and that the fundamental group of $M$ is infinite, does not
contain ${\Z}\oplus {\Z}$, and does not contain a normal cyclic
subgroup, in which case $M$ should be admit a geometric structure
modelled on ${\mathbb H}^3$ [the hyperbolization conjecture].



\subsection{Simplicial volume}
\label{simplicial-volume}

Let $M$ be a closed manifold.  Denote by $C_{*}$ the real chain
complex of $M$: a chain $c\in C_{*}$ is a finite linear combination
$\sum_{i}r_{i}\sigma_{i}$ of singular simplices $\sigma_{i}$ in $M$
with real coefficients $r_{i}$.  Define the {\it simplicial
$l^{1}$-norm} in $C_{*}$ by setting $|c|=\sum_{i}|r_{i}|$. This norm
gives rise to a pseudo-norm on the homology $H_{*}(M,\re)$ by setting
\[| [\alpha] |=\inf\{ |z|\: :\: z\in C_* \mbox{ and } [z] =[\alpha]
\}. \]
When $M$ is orientable, define the {\it simplicial volume} of $M$,
denoted $\|M\|$, to be the simplicial norm of the fundamental class.
The simplicial volume is also called {\it Gromov's invariant}, since
it was first introduced by Gromov \cite{G3}.

The following lower bound on $\| M\|$ is due to Thurston \cite{T}.

\begin{Theorem}\cite[Theorem 6.5.5]{T} Suppose that $M$ is a closed
orientable 3-manifold and that $H\subset M$ is a 3-dimensional
submanifold whose interior admits a complete hyperbolic structure of
finite volume.  Suppose further that $\overline{H}$ is embedded in $M$
and that $\partial\overline{H}$ is incompressible in $M$.  Then,
\[\|M\|\geq \frac{\mbox{\rm Vol}(H)}{v_{3}} >0,\]
where $v_{3}$ is the volume of the regular ideal tetrahedron in
${\mathbb H}^{3}$.
\label{possim}
\end{Theorem}

The next theorem follows immediately from Theorems \ref{possim},
\ref{hypstr}, and \ref{scott1}.

\begin{Theorem} Let $M$ be a closed orientable geometrizable
3-manifold.  Suppose that $\|M\|=0$. Then $M$ is a graph
manifold. 
\label{0graph}
\end{Theorem}

\begin{proof} The proof of Theorem \ref{0graph} is essentially
contained in Soma \cite{So}; we include it here solely for the sake of
completeness.

We begin by considering the prime decomposition of $M$.  That is,
write $M$ as the connected sum $M =M_{0}\#\cdots\# M_{p}$, where each
$M_{i}$ is a prime 3-manifold.   (Note that we are including in this
discussion the case that $M$ is itself prime, and so has trivial prime
decomposition.) 

Since simplicial volume behaves additively with respect to connected
sums (cf. Gromov \cite{G3}), the hypothesis that $M$ has zero
simplicial volume implies that each $M_i$ has zero simplicial volume
as well.  Since the connected sum of graph manifolds is again a graph
manifold (cf. Soma \cite{So}), it suffices to show that each $M_i$ is
a graph manifold.  Since each $M_i$ is prime, it is either irreducible 
or diffeomorphic to ${\mathbb S}^2\times {\mathbb S}^1$, which is a
Seifert fibered space.  So, we may assume without loss of generality
that $M$ is irreducible.

Let $T$ be the torus decomposition of $M$.  Recall that $M$ is assumed
to be geometrizable.  If $T$ is empty, then $M$ admits a geometric
structure other than the one modelled on ${\mathbb H}^3$ (which is
excluded by the assumption on the simplicial volume of $M$), and so
$M$ is a graph manifold, by Theorem \ref{scott1}.

If $T$ is non-empty, then $M$ is sufficiently large, and so Thurston's
geometrization conjecture holds for $M$.  Since $\|M\|=0$, each
component of $M -T$ is a Seifert fibered space, as no piece can be
hyperbolic, by Theorem \ref{possim}.  It follows that $M$ must be a
graph manifold.
\end{proof}

\subsection{Topological entropy}
\label{topological-entropy}

We recall in this subsection the definition of the topological
entropy of the geodesic flow of a smooth Riemannian metric $g$ on a
closed manifold $M$.   For a more detailed discussion, we refer the
interested reader to Paternain \cite{P}.  

The geodesic flow of $g$ is a flow $\phi_{t}$ that acts on $SM$, the
unit sphere bundle of $M$, which is a closed hypersurface of the
tangent bundle of $M$.  Let $d$ be any distance function compatible
with the topology of $SM$. For each $T>0$ we define a new distance
function  
\[d_{T}(x,y):= \max_{0\leq t\leq T}\,d(\phi_{t}(x),\phi_{t}(y)).\]
Since $SM$ is compact, we can consider the minimal number of balls
of radius $\varepsilon>0$ in the metric $d_{T}$ that are necessary
to cover $SM$. Let us denote this number by $N(\varepsilon,T)$. We
define
\[\mbox{\rm h}(\phi,\varepsilon):= \limsup_{T\rightarrow
\infty}\frac{1}{T}\log N(\varepsilon,T).\] 
Observe now that the function $\varepsilon\mapsto \mbox{\rm
h}(\phi,\varepsilon)$ is monotone decreasing and therefore the
following limit exists:
\[\h:= \lim_{\varepsilon\rightarrow 0}\mbox{\rm h}(\phi,\varepsilon).\]
The number $\h$ thus defined is the {\it topological entropy} of the
geodesic flow of $g$. Intuitively, this number measures of orbit
complexity of the flow. The positivity of $h_{top}(\phi)$ indicates
complexity or `chaos' of some kind in the dynamics of
$\phi_{t}$. 

There is a formula, known as Ma\~n\'e's formula, that gives a nice
alternative description of $\h$. Given points $p$ and $q$ in $M$ and
$T>0$, define $n_{T}(p,q)$ to be the number of geodesic arcs joining
$p$ and $q$ with length $\leq T$.  Ma\~n\'e \cite{Man} showed that 
$$\h= \lim_{T\rightarrow \infty}\frac{1}{T}\log \int_{M\times
M}n_{T}(p,q)\;dp\,dq. $$

Finally we note that entropy behaves well under scaling of the metric.
Namely, if $c$ is any positive constant, then ${\rm h}_{{\rm
top}}(cg)=\frac{\h}{\sqrt{c}}$.

\subsection{Minimal volume and collapsing}
\label{minimal-volume}

The {\it minimal volume} MinVol($M$) of a Riemannian manifold $M$ is
defined to be the infimum of $\mbox{\rm Vol}(M,g)$ over all smooth
metrics $g$ such that the sectional curvature $K_{g}$ of $g$ satisfies
$|K_{g}|\leq 1$. This differential invariant was introduced by 
M. Gromov in \cite{G3}.

We shall need the following result, see Cheeger and Gromov
\cite[Example 0.2 and Theorem 3.1]{CG} and Rong \cite{R}.

\begin{Proposition} Let $M$ be a closed orientable 3-manifold.  If $M$
is a graph manifold, then $M$ admits a polarized F-structure, and
hence $\mbox{\rm MinVol}(M)=0$. 
\label{minvol0}
\end{Proposition}

We will not give here the precise definition of a polarized
$F$-structure, because it is too technical.  Instead we give an
informal description, and we refer the interested reader to Cheeger
and Gromov \cite{CG} for a more detailed discussion. 

An {\it $F$-structure} on a manifold $M$ is a natural generalization
of a torus action on $M$. Different tori, possibly of different
dimensions, act on subsets of $M$ in such a way that $M$ is partioned
into disjoint orbits. The $F$-structure is said to be {\it polarized}
if the local actions are locally free. 

Consider the following example of a polarized $F$-structure on a graph
manifold.  Take a compact surface $S$ with non-empty connected
boundary, and consider two copies of $S\times {\mathbb S}^{1}$, each
of which has a 2-torus boundary. Fixing an identification of $\partial
S$ with ${\mathbb S}^1$, glue the boundaries of two copies of $S\times
{\mathbb S}^1$ by a map that interchanges the ${\mathbb S}^{1}$
factors, so that $(x, z)\in \partial S\times {\mathbb S}^1$ on one
copy is glued to $(z,x)\in \partial S\times {\mathbb S}^1$ on the
other copy.

The resulting manifold admits a free circle action on each copy of
${\rm int}(S)\times {\mathbb S}^1$, but at their common boundary the
actions do not agree.  However, they do generate a 2-torus action
which acts locally near their common boundary, thus defining a
polarized $F$-structure on the whole manifold. 

\subsection{An important chain of inequalities}
\label{chain-of-inequalities}

Let $M$ be a closed Riemannian manifold with smooth metric $g$, and
let $\widetilde{M}$ be its universal covering endowed with the induced
metric.  For each $x\in\widetilde{M}$, let $V(x,r)$ be the volume of
the ball with center $x$ and radius $r$.  Set
\[\lambda(g):=\lim_{r\to +\infty}\frac{1}{r}\log\,V(x,r).\]
Manning \cite{Ma} showed that this limit exists and is independent
of $x$.

Set
\[\lambda(M):=\inf\{\lambda(g)\: |\: g \mbox{ is a smooth metric on
$M$ with } {\rm Vol}(M,g) =1\}. \]

It is well known, see Milnor \cite{Mil}, that $\lambda(g)$ is positive
if and only if $\pi_{1}(M)$ has exponential growth. Manning's
inequality \cite{Ma} asserts that for any metric $g$,
\begin{equation}
\lambda(g)\leq {\rm h}_{\rm top}(g).
\label{manineq}
\end{equation}
In particular, it follows that if $\pi_{1}(M)$ has exponential growth,
then ${\rm h}_{\rm top}(g)$ is positive for any metric $g$.  (This
fact was first observed by Dinaburg \cite{D}).  Gromov \cite{G3}
showed that if ${\rm Vol}(M,g) = 1$, then
\begin{equation}
\frac{1}{C_{n}\,n!}\|M\|\leq [\lambda(g)]^n,
\label{simvol-5}
\end{equation}
where
\[C_{n}=\Gamma\left(\frac{n}{2}\right)
\big/\sqrt{\pi}\,\Gamma\left(\frac{n+1}{2}\right).\]
Finally it was observed by Paternain \cite{P} that
\begin{equation}
[{\rm h}(M)]^{n}\leq (n-1)^{n}\mbox{\rm MinVol}(M).
\label{pineq}
\end{equation}
Combining equations (\ref{manineq}), (\ref{simvol-5}), and
(\ref{pineq}), we obtain the following chain of inequalities:
\begin{equation}
\frac{1}{C_{n}n!}\|M\|\leq  [\lambda(M)]^{n}
\leq [{\rm h}(M)]^{n}\leq (n-1)^{n}\mbox{\rm MinVol}(M).
\label{chain-5}
\end{equation}

We note here that the only known 3-manifolds with ${\rm h}(M)>0$
are those with $\|M\|\neq 0$. In fact it follows from Theorem
\ref{0graph}, Proposition \ref{minvol0}, and the chain of
inequalities (\ref{chain-5}) that if $M$ is a closed orientable
geometrizable 3-manifold, then the vanishing of the simplicial volume
implies that ${\rm h}(M)=0$. 

We encode this information in the following theorem.

\begin{Theorem}Let $M$ a closed orientable geometrizable 3-manifold.
Then the following are equivalent:
\begin{enumerate}
\item the minimal volume ${\rm MinVol}(M)$ of $M$ vanishes;
\item the minimal entropy ${\rm h}(M)$ of $M$ vanishes;
\item the simplicial volume $\| M\|$ of $M$ vanishes;
\item $M$ is a graph manifold.
\end{enumerate}
\label{equiv}
\end{Theorem}

\section{Geometric structures and the minimal entropy problem}
\label{geometric-structures}

In this section, we consider the minimal entropy problem for those
3-manifolds that admit a single geometric structure.  Namely, we prove
the following.

\begin{Proposition} Let $M$ be a closed orientable 3-manifold.
Suppose that $M$ admits a geometric structure. Then, the minimal
entropy problem for $M$ can be solved if and only if $M$ admits a
geometric structure modelled on ${\mathbb S}^3$, ${\mathbb E}^3$,
${\mathbb S}^2\times {\re}$, ${\rm Nil}$, or ${\mathbb H}^3$. 
Moreover, if $M$ admits a geometric structure modelled on ${\mathbb
S}^3$, ${\mathbb E}^3$, ${\mathbb S}^2\times {\re}$, or ${\rm Nil}$,
then $M$ admits a smooth metric $g$ with $\h=0$.
\label{preliminary}
\end{Proposition}

\begin{proof}  We start by showing that if $M$ admits a geometric
structure modelled on one of these 5 geometries, then the minimal
entropy problem for $M$ can be solved.  Observe first that if $M$
admits a geometric structure modelled on ${\mathbb H}^{3}$, then the
minimal entropy problem can be solved by the results of Besson,
Courtois and Gallot \cite{BCG}.  

It follows immediately from Theorem \ref{scott1} that if $M$ admits a
geometric stucture modelled on one of the seven geometries ${\mathbb
S}^3$, ${\mathbb E}^3$, ${\mathbb S}^2\times {\re}$, ${\mathbb
H}^{2}\times\re$, $\widetilde{{\rm SL}}_{2}$, ${\rm Nil}$, or ${\rm
Sol}$, then $M$ is a graph manifold.  Hence by Proposition
\ref{minvol0} and the chain of inequalities (\ref{chain-5}), we have
that for such an $M$, the minimal entropy satisfies ${\rm h}(M)=0$.


We now show that if $M$ admits a geometric structure modelled on one
of ${\mathbb S}^3$, ${\mathbb E}^3$, ${\mathbb S}^2\times {\re}$, or 
${\rm Nil}$, then the minimal entropy problem for $M$ can be solved.
To do this, we need to show that $M$ admits a smooth metric $g$ with
$\h =0$.

\begin{enumerate}
\item ${\mathbb S}^3$, ${\mathbb E}^3$, ${\mathbb S}^2\times {\re}$:
All the Jacobi fields in these geometries grow at most linearly (in
the case of ${\mathbb S}^{3}$ they are actually bounded), and hence
all the Liapunov exponents of {\it every} geodesic in $M$ are zero. It
follows from Ruelle's inequality \cite{Ru} that all the measure
entropies are zero.  Hence, by the variational principle, the
topological entropy of the geodesic flow of $M$ must be zero. 
\item ${\rm Nil}$: This geometry can be described as $\re^{3}$ with
the metric
\[ ds^{2}=dx^{2}+dy^{2}+(dz-xdy)^{2}. \]
Here, not all the Jacobi fields grow linearly, but they certainly grow
polynomially. Again this implies that all the Liapunov exponents of
{\it every} geodesic in $M$ are zero and hence the topological entropy
of the geodesic flow of $M$ must be zero.
\end{enumerate}

Since we have assumed that $M$ admits a geometric structure, we
complete the proof by showing that if $M$ admits a geometric structure
modelled on one of remaining geometries, namely ${\mathbb
H}^{2}\times\re$, $\widetilde{{\rm SL}}_{2}$, and ${\rm Sol}$, then
$M$ cannot admit a metric of zero topological entropy.  To do this, we
use the next lemma, together with the fact described in Subsection 
\ref{chain-of-inequalities}, that if $\pi_{1}(M)$ grows exponentially,
then ${\rm h}_{\rm top}(g)>0$ for any smooth metric $g$ on $M$.

\begin{Lemma} Let $M$ be a closed orientable 3-manifold, and suppose
that $M$ admits a geometric structure modelled on one of ${\mathbb
H}^{2}\times\re$, $\widetilde{{\rm SL}}_{2}$, or ${\rm Sol}$.  Then
$\pi_{1}(M)$ grows exponentially.
\label{geom-exp}
\end{Lemma}

\begin{proof} In the case that $M$ admits a geometric structure
modelled on ${\mathbb H}^{2}\times\re$ or $\widetilde{{\rm SL}}_{2}$,
we start by recalling from Theorem \ref{scott1} that $M$ is then a
Seifert fibered space.  The base orbifold of the Seifert fiber space
admits a hyperbolic structure, and so the orbifold fundamental group
of the base orbifold contains a free subgroup of rank 2, and hence so
does $\pi_1(M)$.  Hence, $\pi_1(M)$ grows exponentially. 

In the case that $M$ admits a geometric structure modelled on ${\rm
Sol}$, we have that $M$ is finitely covered by the mapping torus $N$
of a hyperbolic automorphism of a 2-torus.  Note that a hyperbolic
automorphism of a 2-torus is an Anosov diffeomorphism, and so the
suspension flow on $N$ is an Anosov flow. It is known that the
fundamental group of a 3-manifold with an Anosov flow has exponential
growth (see for example Plante and Thurston \cite{PT}).
\end{proof}

This completes the proof of Proposition \ref{preliminary}.
\end{proof}

\section{Proof of Theorem A}
\label{proof-theorem-A}

Up to this point, we have been considering the minimal entropy problem
for closed 3-manifolds that admit a single geometric structure.  In
this section, we consider a more general geometrizable 3-manifold.

\medskip
\noindent
{\bf Theorem A.} {\it Let $M$ be a closed orientable irreducible
3-manifold whose fundamental group contains a ${\Z}\oplus {\Z}$
subgroup. The following are equivalent:
\begin{enumerate}
\item the simplicial volume $\|M\|$ of $M$ is zero and
the minimal entropy problem for $M$ can be solved;
\item $M$ admits a geometric structure modelled on ${\mathbb E}^3$ or ${\rm
Nil}$;
\item $M$ admits a smooth metric $g$ with $\h=0$.
\end{enumerate}}

\begin{proof} Let us show that item 1 implies item 2.  Suppose then
that $M$ has zero simplicial volume and that the minimal entropy
problem for $M$ can be solved.  
We show that $M$ must then admit a geometric structure modelled on
either ${\mathbb E}^3$ or ${\rm Nil}$.  Since the fundamental group of
$M$ contains a ${\Z}\oplus {\Z}$ subgroup, Theorem \ref{link} ensures
that either $M$ contains an incompressible embedded 2-torus or $M$ is
a Seifert fibered space.  We now split the proof into two cases:

\begin{itemize}
\item Suppose first that $M$ contains an incompressible embedded
2-torus, and so is sufficiently large.  Since we have assumed that
$\|M\|=0$, Theorem \ref{0graph} yields that $M$ is a graph manifold.
Hence, by Theorem \ref{equiv}, we have that ${\rm h}(M)=0$.

However, using work of Evans and Moser \cite{evans-moser},
specifically Theorem 4.2 and Corollary 4.10 in \cite{evans-moser}, we
see that either $\pi_1(M)$ contains a free subgroup of rank 2 or $M$
is finitely covered by a 2-torus bundle over ${\mathbb S}^1$.  In the
former case, $\pi_{1}(M)$ grows exponentially and therefore the
minimal entropy problem cannot be solved for $M$.


In the latter case, $M$ admits a geometric structure modelled on one
of ${\mathbb E}^3$, ${\rm Nil}$, or ${\rm Sol}$ (cf. Theorem 5.5 of
Scott \cite{S}).  However, in the case that $M$ admits a geometric
structure modelled on ${\rm Sol}$, we know from Proposition
\ref{preliminary} that the minimal entropy problem cannot be solved
for $M$.

Hence, if the minimal entropy problem can be solved for $M$ and if $M$
contains an incompressible embedded 2-torus, then $M$ admits a
geometric structure modelled on either ${\mathbb E}^3$ or ${\rm Nil}$.

\item The other case is that $M$ is a Seifert fibered space.  Here,
Theorem \ref{scott1} ensures that $M$ possesses a geometric structure
modelled on one of ${\mathbb S}^3$, ${\mathbb E}^3$, ${\mathbb
S}^2\times {\re}$, ${\mathbb H}^{2}\times\re$, $\widetilde{{\rm
SL}}_{2}$ or ${\rm Nil}$.

Since the fundamental group of $M$ admits a ${\Z}\oplus {\Z}$
subgroup, the geometric structure on $M$ cannot be modelled on
${\mathbb S}^3$ or ${\mathbb S}^2\times {\re}$.  Since we have assumed
that the minimal entropy problem can be solved for $M$, Proposition
\ref{preliminary} yields that $M$ must admit a geometric structure
modelled on either ${\mathbb E}^3$ or ${\rm Nil}$, as desired.
\end{itemize}

To see that item 2 implies item 3, recall from Proposition
\ref{preliminary} that if $M$ admits a geometric structure modelled on
${\mathbb E}^3$ or ${\rm Nil}$, then $M$ admits a smooth metric $g$
with $\h=0$. 

Finally to prove that item 3 implies item 1, observe that if $M$
admits a smooth metric $g$ with $\h=0$ it then follows from
inequalities (\ref{manineq}) and (\ref{simvol-5}) that $M$ has zero
simplicial volume. 

This completes the proof of Theorem A.
\end{proof}

\section{Proof of Theorem B}
\label{proof-theorem-B}

We are now ready to consider the minimal entropy problem for a general
geometrizable 3-manifold with zero simplicial volume.

\medskip
\noindent
{\bf Theorem B.} {\it Let $M$ be a closed orientable geometrizable
3-manifold. The following are equivalent:
\begin{enumerate}
\item the simplicial volume $\|M\|$ of $M$ is zero and
the minimal entropy problem for $M$ can be solved;
\item $M$ admits a geometric structure modelled on  ${\mathbb S}^{3}$, ${\mathbb
S}^{2}\times \re$, ${\mathbb E}^3$, or ${\rm Nil}$;
\item $M$ admits a smooth metric $g$ with $\h=0$.
\end{enumerate}}

\begin{proof} Let us prove that item 1 implies item 2.  Suppose that
$M$ has zero simplicial volume and that the minimal entropy problem
for $M$ can be solved. Since $M$ is geometrizable and its simplicial
volume vanishes, Theorem \ref{0graph} tell us that $M$ is a graph
manifold. Hence, by Theorem \ref{equiv}, $M$ has zero minimal
entropy.

Since we are assuming that the minimal entropy problem can be solved
for $M$, the fact that $M$ has zero minimal entropy in turn implies
there exists a smooth metric on $M$ with zero topological entropy.
This in turn implies, by the discussion in Section
\ref{chain-of-inequalities}, that $\pi_{1}(M)$ does not have
exponential growth.

However, it is a fact from combinatorial group theory (which
follows immediately from the existence of normal forms for free
products, for instance) that if $A$ and $B$ are two finitely generated
groups, then the free product $A*B$ contains a free subgroup of rank
two unless $A$ is trivial or $B$ is trivial, or $A$ and $B$ are both
of order two. Since the fundamental group of a connected sum is 
the free product of the fundamental groups of the summands, we
conclude that either the prime decomposition is trivial or there are
only two summands both of which have fundamental group $\Z_{2}$.

In the former case, it follows that $M$ must be either irreducible or
${\mathbb S}^{2}\times {\mathbb S}^{1}$, while in the latter case $M$
must be ${\mathbb P}^{3}\# {\mathbb P}^{3}$, where ${\mathbb P}^3$ is
the 3-dimensional real projective space.  Since ${\mathbb S}^{2}\times
{\mathbb S}^{1}$ and ${\mathbb P}^{3}\# {\mathbb P}^{3}$ both admit a
geometric structure modelled on ${\mathbb S}^{2}\times \re$, we may
assume from now on that $M$ is irreducible.

There are now several cases, depending on $\pi_1(M)$.  Suppose first
that $\pi_1(M)$ is finite.  Since $M$ is geometrizable, we have that
$M$ admits a geometric structure modelled on ${\mathbb S}^3$.

In the case that $\pi_1(M)$ is infinite and contains a ${\mathbb
Z}\oplus {\mathbb Z}$ subgroup, the assumption that the simplicial
volume of $M$ is zero, together with the fact that the minimal entropy
problem can be solved for $M$, allows us to apply Theorem A to see
that $M$ admits a geometric structure modelled on ${\mathbb E}^3$ or
${\rm Nil}$.

The remaining case is that $\pi_1(M)$ is infinite and does not contain
a ${\mathbb Z}\oplus {\mathbb Z}$ subgroup.  Since $M$ is
geometrizable, either $M$ admits a hyperbolic structure or $M$ is
Seifert fibered.  (Since $\pi_1(M)$ does not contain a ${\mathbb
Z}\oplus {\mathbb Z}$ subgroup, $M$ cannot admit a geometric structure
modelled on ${\rm Sol}$, as ${\rm Sol}$ manifolds are finitely covered
by 2-torus bundles over the circle.)  However, since $\| M\| =0$, $M$
cannot admit a hyperbolic structure.

Note though that $M$ cannot admit a geometric structure modelled on
${\mathbb H}^2\times {\mathbb R}$, ${\mathbb E}^3$, $\widetilde{{\rm
SL}_2}$, or ${\rm Nil}$, as such manifolds always have a ${\mathbb
Z}\oplus {\mathbb Z}$ in their fundamental groups, by Lemma
\ref{seifert-fundamental-group}.  Hence, the only possibilities
remaining are that $M$ admits a geometric structure modelled on either
${\mathbb S}^2\times {\mathbb R}$ or ${\mathbb S}^3$, as desired.

To see that item 2 implies item 3, recall from Proposition
\ref{preliminary} that if $M$ admits a geometric structure modelled on
${\mathbb S}^{3}$, ${\mathbb S}^{2}\times \re$, ${\mathbb E}^3$, or
${\rm Nil}$, then $M$ admits a smooth metric $g$ with $\h=0$. 

Finally to prove that item 3 implies item 1, observe that if $M$
admits a smooth metric $g$ with $\h=0$, it then follows from
inequalities (\ref{manineq}) and (\ref{simvol-5}) that $M$ has zero
simplicial volume.

This completes the proof of Theorem B.
\end{proof}

\end{document}